\documentclass[conference]{IEEEtran}

\ifCLASSINFOpdf
\else

\fi

\usepackage[cmex10]{amsmath}
\usepackage{amssymb}

\usepackage{graphics} 
\usepackage{epsfig}

\usepackage{eso-pic}

\begin{document}
%
\title{General formula for stability testing of fractional-delay systems}


\author{\IEEEauthorblockN{Farshad Merrikh-Bayat}\\
\IEEEauthorblockA{Faculty of Electrical Engineering\\
University of Zanjan\\
Zanjan, Iran, P.O.Box 313\\
Email: f.bayat@znu.ac.ir}

}

\maketitle


\begin{abstract}
An easy-to-use and effective formula for stability testing of a
system with fractional-delay characteristic equation in the
general form of $\Delta(s)=P_0(s)+\sum_{i=1}^N P_i(s)\exp(-\zeta_i
s^{\beta_i}) =0$, where $P_i(s)$ ($i=0,\ldots, N$) are the
so-called fractional-order polynomials and $\zeta_i$ and $\beta_i$
are positive real constants, is proposed in this paper. The
proposed formula determines the number of unstable roots of the
characteristic equation (i.e., those located in the right
half-plane of the first Riemann sheet) by applying Rouche's
theorem. Numerical simulations are also presented to confirm the
efficiency of the proposed formula.
\end{abstract}


%
\IEEEpeerreviewmaketitle

\section{Introduction}
In some of the recently-developed control problems we need to
check the stability of a system with the so-called
\emph{fractional-delay} characteristic equation in the general
form of
\begin{equation}\label{frac_delay}
\Delta(s)=P_0(s)+\sum_{i=1}^N P_i(s)\exp(-\zeta_i s^{\beta_i}) =0,
\end{equation}
where $\zeta_i$ and $\beta_i$ are positive real constants and
$P_i(s)$ ($i=1,\ldots,N$) are fractional-order polynomials in the
form of
\begin{equation}
P_i(s)=\sum_{k=1}^{M_i}a_{ik}s^{\alpha_{ik}},
\end{equation}
where $a_{ik}$ and $\alpha_{ik}$ are real and positive real
constants, respectively, and
\begin{equation}\label{char}
P_0(s)=s^{\alpha_n}+a_{n-1}s^{\alpha_{n-1}}+\ldots+a_1s^{\alpha_1}+a_0,
\end{equation}
where, without any loss of generality, it is assumed that the
powers of $s$ in (\ref{char}) satisfy the following relations:
\begin{equation}
\alpha_n>\alpha_{n-1}>\dots>\alpha_1>0.
\end{equation}
As an example of a system with fractional-delay characteristic
equation, consider a classical unity-feedback system in which a
process with transfer function
\begin{equation}
G(s)=\frac{K}{1+sT}e^{-sL},
\end{equation}
is controlled with the so-called $\mathrm{PI}^\lambda
\mathrm{D}^\mu$ controller with transfer function
\cite{podlubny99a}:
\begin{equation}
C(s)=K_p\left(1+\frac{1}{T_is^\lambda}+T_d s^\mu\right),
\end{equation}
where $K_P$, $T_i$, $T_d$, $\lambda$, and $\mu$ are unknown
parameters of the controller to be determined. It can be easily
verified that the characteristic equation of this system is as the
following
\begin{equation}
\Delta(s)=T_is^\lambda (1+sT)+K_pK(T_is^\lambda +1+T_iT_d
s^{\lambda+\mu})e^{-sL}=0,
\end{equation}
which can be considered as a special case of (\ref{frac_delay})
with $N=1$,
\begin{equation}
P_0(s)=T_is^\lambda (1+sT),
\end{equation}
\begin{equation}
P_1(s)=K_pK(T_is^\lambda +1+T_iT_d s^{\lambda+\mu}),
\end{equation}
$\zeta_1=L$, and $\beta_1=1$. If in this example one tries to find
the optimal values of $K_P$, $T_i$, $T_d$, $\lambda$, and $\mu$ by
means of meta-heuristic optimization algorithms such that a
certain cost function (e.g., ISE performance index corresponding
to the tracking of unit step command) is minimized, he/she will
need a method to check the feasibility of the solutions generated
by the meta-heuristic optimization algorithm from the stability
point of view. It should be noted that since in such optimization
problems the cost function is usually expressed in the frequency
domain (by applying Parseval's theorem), the resulted optimal
controller may destabilize the feedback system \cite{farshad3}.

As a more complicated example, consider the problem of designing
an optimal $\mathrm{PI}^\lambda \mathrm{D}^\mu$ controller for a
process whose transfer function consists of fractional powers of
$s$ possibly in combination with exponentials of fractional powers
of $s$. For example, the transfer functions:
\begin{equation}\label{fd2}
G(s)=\frac{\cosh\left(x_0\sqrt{s}\right)}{\sqrt{s}\sinh \left(
\sqrt{s}\right)},\quad 0<x_0<1,
\end{equation}
and
\begin{equation}\label{fd3}
G(s)=\frac{\sinh\left(x_0\sqrt{s}\right)}{\sinh \left(
\sqrt{s}\right)},\quad 0<x_0<1,
\end{equation}
appear in boundary control of one-dimensional heat equation with
Neumann and Dirichlet boundary conditions \cite{zwart}. Other
examples of this type can be found in \cite{zwart}-\cite{helie}.
Moreover, in some applications in order to arrive at more accurate
models, the process is modelled with a fractional-order transfer
function. For instance, Podlubny \cite{podlubny99} showed that the
fractional-order transfer function:
\begin{equation}\label{fd4}
G(s)=\frac{1}{0.7943s^{2.5708}+5.2385s^{0.8372}+1.5560},
\end{equation}
can better model a heating furnace compared to classical
integer-order transfer functions. Clearly, in dealing with
complicated transfer functions such as those given in
(\ref{fd2})-(\ref{fd4}) we need more powerful tools to determine
the stability of the corresponding closed-loop system.

Stability analysis of the feedback system when such complicated
transfer functions exist in the loop is a challenging task. Even
the stability analysis of a feedback system which consists of both
$\mathrm{PI}^\lambda \mathrm{D}^\mu$ controller and a process with
dead-time is not straightforward. So far, many researchers have
tried to develop analytical or numerical methods for stability
testing of systems with fractional-delay characteristic equations
(see \cite{hwang} for a detailed review of some important works in
relation to the stability testing of fractional-delay systems).
Probably, the most famous analytical method for stability testing
of fractional-order systems (as a special case of fractional-delay
systems) is the \emph{sector stability test} of Matignon
\cite{matignon1}, which was already reported in the work of Ikeda
and Takahashi \cite{ikeda}. Application of this method is limited
to the case where the sigma term does not exist in
(\ref{frac_delay}) and $P_0(s)$ is of \emph{commensurate order}.
Few numerical algorithms for stability testing of
(\ref{frac_delay}) can also be found in the literature (see, for
example, \cite{farshad1} and \cite{hwang} and the references
therein for more information on this subject). As far as we know,
all of these methods suffer from the limitation that can be
applied only to a certain class of fractional-delay systems
\cite{farshad1}, or the results are of probabilistic nature
\cite{hwang}.

The aim of this paper is to propose a formula for determining the
number of unstable roots of (\ref{frac_delay}). The proposed
formula is actually a generalization of the method already
proposed by author in \cite{farshad1}. However, the formula
developed in this paper has the advantage of being much simpler
compared to the one presented in \cite{farshad1}, and moreover, it
can be easily applied to a more general form of fractional-delay
characteristic equations.

The rest of this paper is organized as follows. The proposed
formula for stability testing of fractional-delay characteristic
equations is presented in Section \ref{sec_propos}. Four numerical
examples are studied in Section \ref{sec_exam}, and finally,
Section \ref{sec_conc} concludes the paper.

\section{Proposed formula for stability testing of fractional-delay characteristic
equations}\label{sec_propos} The first step in dealing with
multi-valued complex functions (such as the one presented in
(\ref{frac_delay})) is to construct the domain of definition of
the function appropriately. The domain of definition of the
characteristic function given in (\ref{frac_delay}) is, in
general, in the form of a Riemann surface with infinite number of
Riemann sheets, where the origin is a branch point and the branch
cut is considered (arbitrarily) at $\mathbb{R}^-$. Equation
$\Delta(s)=0$ as defined in (\ref{frac_delay}) has, in general,
infinite number of roots which are distributed on this Riemann
surface. As a well-known fact, a system with characteristic
equation (\ref{frac_delay}) is stable if and only if it does not
have any roots in the right half-plane of the first Riemann sheet
\cite{hwang,farshad2}. Hence, stability analysis of a system with
characteristic equation (\ref{frac_delay}) is equivalent to
investigation for the roots of $\Delta(s)=0$ in the right
half-plane of the first Riemann sheet. In the following we will
use Rouche's theorem for this purpose.

First, let us briefly review the Rouche's theorem. Consider the
complex function $f:\mathbb{C}\rightarrow\mathbb{C}$ which has
zeros of orders $m_1,\ldots,m_k$ respectively at $z_1,\ldots,z_k$
and does not have any poles. This function can be written as
\begin{equation}\label{f}
f(s)=g(s)(s-z_1)^{m_1}\times (s-z_2)^{m_2}\times \cdots \times
(s-z_k)^{m_k},
\end{equation}
where $g(s)$ has neither pole nor zero. Taking the natural
logarithm from both sides of (\ref{f}) leads to
\begin{multline}\label{lnf}
\ln f(s)=\ln g(s)+m_1 \ln (s-z_1)+m_2 \ln (s-z_2)+\ldots\\
+m_k\ln(s-z_k).
\end{multline}
Derivation with respect to $s$ from both sides of (\ref{lnf})
yields
\begin{equation}
\frac{f'(s)}{f(s)}=\frac{g'(s)}{g(s)}+\frac{m_1}{s-z_1}+\frac{m_2}{s-z_2}
+\ldots+\frac{m_k}{s-z_k}.
\end{equation}
Now let $\gamma$ be a simple, closed, counterclockwise contour
such that $f(s)$ has no zeros (or singularities like branch point
and branch cut in dealing with multi-valued functions) on it. Then
it is concluded from the Residue theorem that
\begin{equation}\label{rouch1}
\frac{1}{2\pi i}\oint_\gamma\frac{f'(s)}{f(s)}\mathrm{d}s=M,
\end{equation}
where $M$ is equal to the total number of the roots of $f(s)=0$
inside $\gamma$. Clearly, if the contour $\gamma$ is considered
such that all zeros of $f(s)$ lie inside it then we have
$M=\sum_{j=1}^k m_j$. Equation (\ref{rouch1}) can be used to
calculate the number of zeros of the given function $f(s)$ inside
the desired contour $\gamma$ (which, of course, should have the
above-mentioned properties). For this purpose, we can simply use a
numerical integration technique to evaluate the integral in the
right hand side of (\ref{rouch1}) for the given contour $\gamma$
and function $f$.

\begin{figure}[tb]
\begin{center}
\includegraphics[width=5.5cm]{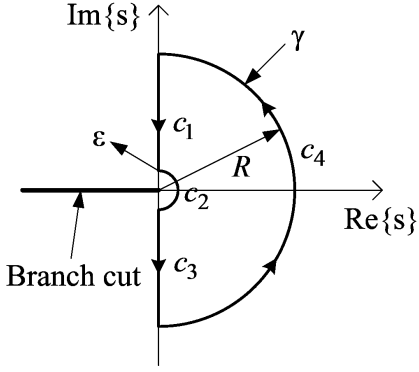}
\caption{The contour $\gamma$ considered on the first Riemann
sheet for stability testing of the fractional-delay system under
consideration. A system with characteristic equation $\Delta(s)=0$
(as defined in (\ref{frac_delay})) is stable if and only if it
does not have any roots inside $\gamma$.} \label{fig_contour}
\end{center}
\end{figure}

According to the above discussions, by setting $f(s)$ equal to
$\Delta(s)$ and $\gamma$ equal to the border of the region of
instability (which is equal to the closed right half-plane of the
first Riemann sheet) the value obtained for $M$ from
(\ref{rouch1}) will be equal to the number of unstable roots of
the characteristic equation. In the following, we consider the
contour $\gamma$ as shown in Fig. \ref{fig_contour} and
$f(s)=\Delta(s)$ (where $\Delta(s)$ is defined in
(\ref{frac_delay})) and then simplify the integral in the left
hand side of (\ref{rouch1}) to arrive at a more effective formula
for stability testing of the fractional-delay system under
consideration (clearly, the system is stable if and only if
$M=0$). Note that the very small semicircle in Fig.
\ref{fig_contour} is used to avoid the branch-point located at the
origin.

According to (\ref{rouch1}) and Fig. \ref{fig_contour} we can
write
\begin{equation}\label{m}
M=\frac{1}{2\pi
i}\oint_{\gamma}\frac{\Delta'(s)}{\Delta(s)}\mathrm{d}s =
\frac{1}{2\pi i}\left(\int_{c_1+c_3} + \int_{c_2}+ \int_{c_4}
\right).
\end{equation}
In (\ref{m}), the integral $\int_{c_1+c_3}$ is calculated as
\begin{align}
\int_{c_1+c_3}&=\int^\varepsilon_\infty
\frac{\Delta'(i\omega)}{\Delta(i\omega)} i\mathrm{d}\omega +
\int_{-\varepsilon}^{-\infty}
\frac{\Delta'(i\omega)}{\Delta(i\omega)} i\mathrm{d}\omega\\
&=-\int_\varepsilon^\infty
\frac{\Delta'(i\omega)}{\Delta(i\omega)} i\mathrm{d}\omega+
\int_\varepsilon^\infty \frac{\Delta'(-i\omega)}{\Delta(-i\omega)}
(-i)\mathrm{d}\omega\\ &=-i\int_\varepsilon^\infty
\frac{\Delta'(i\omega)}{\Delta(i\omega)} \mathrm{d}\omega-i
\int_\varepsilon^\infty
\left(\frac{\Delta'(i\omega)}{\Delta(i\omega)} \right)^\ast
\mathrm{d}\omega,
\end{align}
which yields
\begin{equation}\label{c1}
\int_{c_1+c_3}=-2i\int_\varepsilon^\infty
\mathrm{Re}\left\{\frac{\Delta'(i\omega)}{\Delta(i\omega)}
\right\} \mathrm{d}\omega.
\end{equation}
The integral $\int_{c_2}$ in (\ref{m}) is calculated as
\begin{align}\label{c2}
\int_{c_2}&=\lim_{\varepsilon\rightarrow 0}\int_{\frac {\pi}
{2}}^{-\frac{\pi}{2}}\frac{\Delta'(\varepsilon e^{i\theta})}
{\Delta(\varepsilon e^{i\theta})}\varepsilon ie^{i\theta}
\mathrm{d}\theta\\&= \int_{\frac {\pi} {2}}^{-\frac{\pi}{2}}
\lim_{\varepsilon\rightarrow 0} \left\{\varepsilon
\frac{\Delta'(\varepsilon e^{i\theta})} {\Delta(\varepsilon
e^{i\theta})}\right\} ie^{i\theta} \mathrm{d}\theta.
\end{align}
In the above equation $\lim_{\varepsilon\rightarrow 0}
\Delta(\varepsilon e^{i\theta})$ is equal to a nonzero constant
(else, the characteristic function has a strong singularity at the
origin and the corresponding system is unstable) and $\Delta'
(\varepsilon e^{i\theta})\sim K\varepsilon^\eta$ as
$\varepsilon\rightarrow 0$, where $K$ and $\eta>-1$ are two
constants. Hence, $\int_{c_2}$ tends to zero as $\varepsilon
\rightarrow 0$. Finally, $\int_{c_4}$ in (\ref{m}) is calculated
as
\begin{align}
\int_{c_4}&=\lim_{R\rightarrow\infty}\int_{-\frac {\pi}
{2}}^{\frac{\pi}{2}}\frac{\Delta'(Re^{i\theta})}
{\Delta(Re^{i\theta})}Rie^{i\theta} \mathrm{d}\theta\\&
=\int_{-\frac {\pi} {2}}^{\frac{\pi}{2}}
\lim_{R\rightarrow\infty}\left\{\frac{\Delta'(Re^{i\theta})}
{\Delta(Re^{i\theta})}R\right\} ie^{i\theta} \mathrm{d}\theta\\ &=
i\alpha_n \pi.\label{c4}
\end{align}
(See (\ref{frac_delay}) and (\ref{char}) for the definition of
$\alpha_n$.) Substitution of (\ref{c1}) and (\ref{c4}) in
(\ref{m}) and considering the fact that $\int_{c_2}=0$ results in
\begin{equation}\label{main}
M=\frac{\alpha_n}{2}-\frac{1}{\pi}\int_{\varepsilon=0^+}^\infty
\mathrm{Re}\left\{\frac{\Delta'(i\omega)}{\Delta(i\omega)}
\right\} \mathrm{d}\omega,
\end{equation}
where $M$ is equal to the number of unstable poles of a system
with characteristic equation $\Delta(s)=0$ as defined in
(\ref{frac_delay}).

Equation (\ref{main}) is the main result of this paper. It should
be noted that the value of $\varepsilon$ in (\ref{main}) cannot,
in general, be considered exactly equal to zero. That is because
of the fact that the numerical integration technique used to
evaluate the integral in (\ref{main}) performs this task by
evaluating the integrand at different points of the $\omega$ axis.
Hence, the numerical integration algorithm may halt if the
integrand becomes singular at the origin (which is the case if,
for example, $0<\alpha_1<1$ in (\ref{char})). In practice, in
order to determine the number of unstable poles of the given
fractional-delay transfer function we can consider the lower and
upper bound of the integral in (\ref{main}) equal to sufficiently
small and big positive numbers, respectively. The MATLAB function
\textsf{quadl} (as well as \textsf{quadgk}) can be used to
evaluate the integral in (\ref{main}). Some numerical examples
will be presented in the next section.

\section{Numerical examples}\label{sec_exam}
In the following we study the application of (\ref{main}) for
stability testing of some fractional-delay systems. In each case,
the impulse response of the corresponding system is also plotted
to verify the correctness of the result. The method used in this
paper to calculate the impulse response of the given
fractional-order system is based on the formula proposed in
\cite{valsa} for numerical inversion of Laplace transforms. In
this method the impulse response of the given fractional-order
system is approximated by numerical inversion of its transfer
function. The MATLAB code of this method, \textsf{invlap.m}, can
freely be downloaded from
http://www.mathworks.com/matlabcentral/fileexchange/. Most of the
following examples have already been studied by author in
\cite{farshad1}.

\textbf{Example 1.} Consider a system with characteristic equation
\begin{align}
\Delta_1(s)&=(s^{\pi/2}+1)(s^{\pi/3}+1)\\
&=s^{5\pi/6}+s^{\pi/2}+s^{\pi/3}+1=0.\label{char1}
\end{align}
The roots of this equation can be calculated analytically, which
are as the following:
\begin{equation}\label{r1}
s_{k_1}=e^{j2(2k_1+1)}, \quad k_1\in\mathbb{Z},
\end{equation}
and
\begin{equation}\label{r2}
s_{k_2}=e^{j3(2k_2+1)},\quad k_2\in\mathbb{Z}.
\end{equation}
As it is observed, the characteristic equation given in
(\ref{char1}) has infinite many roots which are distributed on a
Riemann surface with infinite number of Riemann sheets. It is
concluded from (\ref{r1}) and (\ref{r2}) that (\ref{char1}) has
four roots on the first Riemann sheet which are $e^{\pm j 2}$ and
$e^{\pm j3}$, and none of them are located in the right half-plane
(recall that all roots whose phase angle lies in the range
$[-\pi,\pi)$ belong to the first Riemann sheet).

\begin{figure}[tb]
\begin{center}
\includegraphics[width=8.5cm]{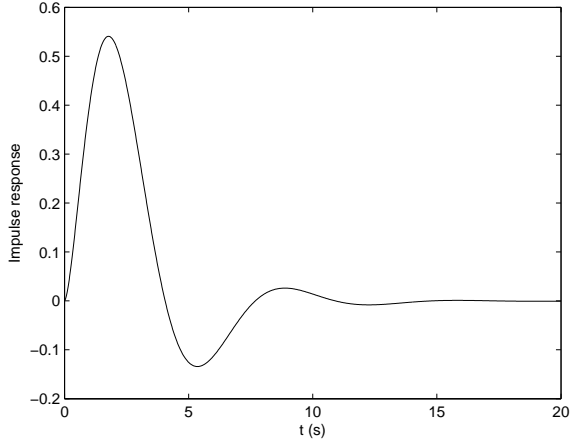}
\caption{Impulse response of a system with transfer function
(\ref{h1}).} \label{fig_ex1}
\end{center}
\end{figure}

Comparing (\ref{char1}) with (\ref{frac_delay}) and (\ref{char})
yields $\alpha_n=5\pi/6$ (note that (\ref{char1}) has no delay
terms). Application of (\ref{main}) assuming that the lower and
upper bound of integral in (\ref{main}) are equal to 0 and 1000,
respectively, yields $M=2.3300\times 10^{-4}$ which is consistent
with the above-mentioned analytical result. Figure \ref{fig_ex1}
shows the impulse response of a system with transfer function
\begin{equation}\label{h1}
H_1(s)=\frac{1}{\Delta_1(s)}=\frac{1}
{s^{5\pi/6}+s^{\pi/2}+s^{\pi/3}+1}.
\end{equation}
As it can be observed in this figure, the impulse response of the
system is absolutely summable, as it is expected.\\

\begin{figure}[tb]
\begin{center}
\includegraphics[width=8.5cm]{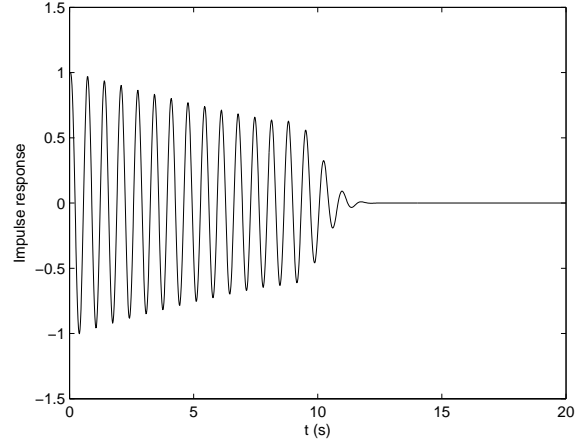}
\caption{Impulse response of a system with transfer function
(\ref{h2}).} \label{fig_ex2}
\end{center}
\end{figure}

\textbf{Example 2.} Stability of a system with fractional-delay
characteristic equation:
\begin{equation}
\Delta_2(s)=s+K(\sqrt{s}+1)e^{-\sqrt{s}}=0,
\end{equation}
is studied in \cite{ozturk} and it is especially shown that it is
stable for $K<21.51$ and unstable for $K>21.51$. Application of
(\ref{main}) assuming that $K=21$, $\alpha_n=1$, and the lower and
upper bound of integral are equal to 0 and 500, respectively,
yields $M=3.4227\times 10^{-9}$, which implies the stability of
system as it is expected. Figure \ref{fig_ex2} shows the impulse
response of a system with transfer function
\begin{equation}\label{h2}
H_2(s)=\frac{1}{\Delta_2(s)}=\frac{1}{s+21(\sqrt{s}+1)e^{-\sqrt{s}}}.
\end{equation}
As it can be observed, the impulse response is absolutely
summable, as it is expected.
Repeating the above procedure with $K=22$ yields $M=2.0174$, which
means that in this case the system has two unstable poles. This
result is also consistent with the one presented in
\cite{farshad1}.\\

\begin{figure}[tb]
\begin{center}
\includegraphics[width=8.5cm]{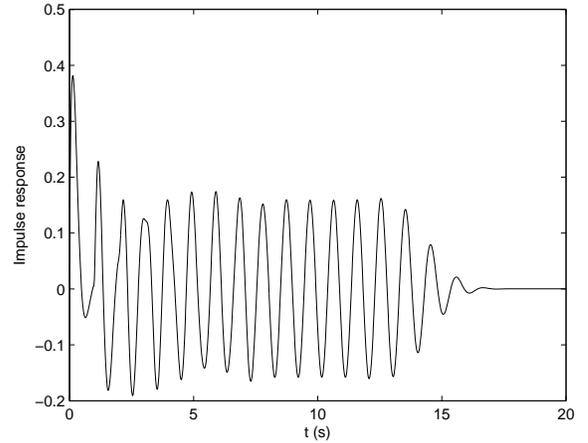}
\caption{Impulse response of a system with transfer function
(\ref{h3}).} \label{fig_ex3}
\end{center}
\end{figure}

\textbf{Example 3.} It is shown in \cite{ozturk2} that a system
with characteristic equation
\begin{equation}
\Delta_3(s)=s^{1.5}-1.5s+4s^{0.5}+8-1.5s e^{-\tau s}=0,
\end{equation}
is stable for the values of $\tau\in(0.99830, 1.57079)$ and
unstable for other values of $\tau$. It is also shown by author in
\cite{farshad1} that this system has two unstable poles for
$\tau=0.99$. Application of (\ref{main}) assuming $\tau=1$ and
considering the fact that here we have $\alpha_n=1.5$ yields
$M=0.0082$ (the lower and upper bound of integral are considered
equal to 0 and 500, respectively). As it is observed, the result
obtained by using the proposed method is fairly close to zero.
Figure \ref{fig_ex3} shows the impulse response of a system with
transfer function:
\begin{equation}\label{h3}
H_3(s)=\frac{1}{\Delta_3(s)}=\frac{1}{s^{1.5}-1.5s+4s^{0.5}+8-1.5s
e^{-s}}.
\end{equation}
As it can be observed in this figure, the impulse response of the
system is absolutely summable and consequently, the corresponding
system is stable.
In this example, application of (\ref{main}) assuming $\tau=0.99$
yields $1.9994$ which is consistent with the result presented in
\cite{farshad1}.\\

\begin{figure}[tb]
\begin{center}
\includegraphics[width=8.5cm]{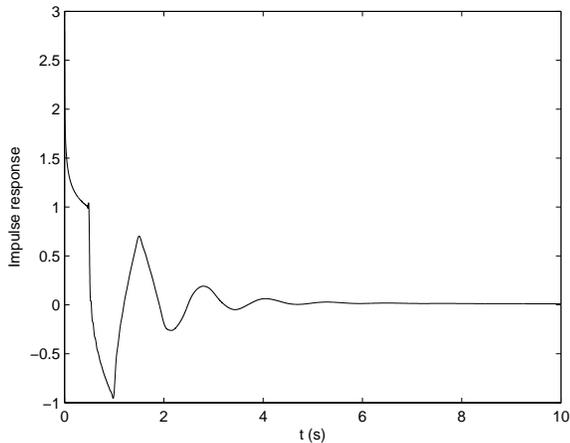}
\caption{Impulse response of a system with transfer function
(\ref{h4}).} \label{fig_ex4}
\end{center}
\end{figure}

\textbf{Example 4.} It is shown in \cite{hwang} (by applying
Lambert W function) that a system with the following
characteristic equation
\begin{equation}
\Delta_4(s)=s^{5/6}+(s^{1/2}+s^{1/3})e^{-0.5s}+e^{-s}=0,
\end{equation}
is stable. Clearly, here we have $\alpha_n=5/6$. Application of
(\ref{main}) (assuming that the lower and upper bound of integral
are equal to 0 and 100, respectively) leads to $M=0.0290$, which
implies the stability of system. Figure \ref{fig_ex4} shows the
impulse response of a system with transfer function
\begin{equation}\label{h4}
H_4(s)=\frac{1}{\Delta_4(s)}=\frac{1}{s^{5/6}+(s^{1/2}+s^{1/3})e^{-0.5s}+e^{-s}}.
\end{equation}
It can be observed that the impulse response is absolutely
summable and consequently, the system is stable, as it is
expected.

\section{Conclusion}\label{sec_conc}
An easy-to-use, effective and very general formula for stability
testing of fractional-delay systems is proposed in this paper. The
proposed formula can be used to determine the number of unstable
poles of a system whose characteristic equation contains, in
general, both fractional powers of $s$ and exponentials of
fractional powers of $s$.

\end{document}